\documentclass{sig-alternate}

\newtheorem{definition}{Definition}
\newtheorem{theorem}{Theorem}
\newtheorem{problem}{Problem} 
\newtheorem{remark}{Remark}
\newtheorem{example}{Example}

\newcommand{\NN}{\mathbb{N}}
\newcommand{\QQ}{\mathbb{Q}}
\newcommand{\ZZ}{\mathbb{Z}}

\newcommand{\KK}{\ensuremath{\mathbb{K}}}
\newcommand{\GG}{\ensuremath{\mathcal{G}}}
\newcommand{\FF}{\ensuremath{\mathcal{F}}}
\newcommand{\lm}{\mathrm{lm}}
\newcommand{\lc}{\mathrm{lc}}
\newcommand{\Hom}{\mathrm{Hom}}
\newcommand{\Sol}{\mathrm{Sol}}
\newcommand{\End}{\mathrm{End}}

\usepackage{algorithm}
\usepackage{algorithmic}
\usepackage{url}

\begin{document}
%
\conferenceinfo{ISSAC}{'16 Waterloo, Ontario, Canada}

\title{A Factorization Algorithm for
  {\fontsize{1.05em}{1em}{$G$}}-Algebras and Applications}
%
%
%
%
%

\numberofauthors{2} 
%
\author{
%
%
\alignauthor
Albert Heinle\\
       \affaddr{David R. Cheriton School of Computer Science}\\
       \affaddr{University of Waterloo}\\
       \affaddr{200 University Avenue West}\\
       \affaddr{Waterloo, ON N2L 3G1, Canada}\\
       \email{aheinle@uwaterloo.ca}
\alignauthor
Viktor Levandovskyy\\
       \affaddr{Lehrstuhl D f\"ur Mathematik}\\
       \affaddr{RWTH Aachen University}\\
       \affaddr{Pontdriesch 16}\\
       \affaddr{52062 Aachen, Germany}\\
       \email{levandov@math.rwth-aachen.de}
}

\maketitle
\begin{abstract}
It has been recently discovered by Bell, Heinle and Levandovskyy that
a large class of algebras, including the ubiquitous $G$-algebras,
are finite factorization domains (FFD for short).

Utilizing this result, we contribute an algorithm to find all distinct
factorizations of a given element $f \in \GG$, where $\GG$ is any
$G$-algebra, with minor assumptions on the underlying field.

Moreover, the property of being an FFD, in combination with the
factorization algorithm, enables us to propose an analogous description
of the factorized Gr\"obner basis algorithm for $G$-algebras. This
algorithm is useful for various applications, e.g.  in analysis of
solution spaces of systems of linear partial functional equations with
polynomial coefficients, coming from $\GG$.  Additionally, it is
possible to include inequality constraints for ideals in the input.
\end{abstract}




\section{Introduction}

\emph{Notations:} Throughout the paper we denote by $\KK$ a field. In
the algorithmic part we will assume $\KK$ to be a computable
field. $\NN_0 = \NN \cup \{0\}$ is the set of natural numbers
including zero.
For a $\KK$-algebra $R$
we denote by $U(R)$ the group of invertible (unit) elements of $R$, which is
nonabelian in general. For $f\in R$ we denote by $Rf$ the left ideal, generated by $f$. 
The main focus in this paper lies in so called
$G$-algebras, which are defined as follows.

\begin{definition}
\label{def: G-algebra}
For $n\in\NN$ and $1\leq i < j \leq n$ consider the units $c_{ij} \in \KK^*$ and 
polynomials $d_{ij} \in \KK[x_1,\ldots,x_n]$. Suppose, that there exists a
monomial total  well-ordering $\prec$ on $\KK[x_1,\ldots,x_n]$, such that 
for any $1\leq i < j \leq n$ either $d_{ij}=0$ or the leading monomial of $d_{ij}$ is smaller than $x_i x_j$ with respect to $\prec$. The $\KK$-algebra
$A := \KK\langle x_1,\ldots, x_n \mid \{ x_j x_i = c_{ij} x_i x_j + d_{ij} \colon  1\leq i < j \leq n \} \rangle$ is called a $G$-\emph{algebra}, if 
$\{x_1^{\alpha_1} \cdot \ldots \cdot x_n^{\alpha_n} : \alpha_i\in\NN_0 \}$ is a $\KK$-basis of $A$.
\end{definition}
$G$-algebras \cite{AP,LVDiss} are also known as algebras of solvable type \cite{KW,HLi} and as PBW algebras \cite{GomezRefiltering,Bueso:2003}.
$G$-algebras are Noetherian domains of finite global, Krull and Gel'fand-Kirillov
dimensions.

We assume that the reader is familiar with the basic terminology in
the area of Gr\"obner bases, both in the commutative as well as in
the non-commutative case. We recommend
\cite{Buchberger:1997,Bueso:2003,LVDiss} as literature on this topic.

Recall, that $r\in R\setminus\{0\}$ is called {\bf irreducible}, if in any factorization
$r = ab$ either $a\in U(R)$ or $b\in U(R)$ holds. Otherwise, we call 
$r$ reducible.

\begin{definition}[cf. \cite{bell2014noncommutative}]
\label{def: FFD}
  Let $A$ be a (not necessarily commutative) domain.  We say that $A$ is a
  \textit{finite factorization domain} (FFD, for short), if
  every nonzero, non-unit element of $A$ has at least one
  factorization into irreducible elements and there are at most
  finitely many distinct factorizations into irreducible elements up
  to multiplication of the irreducible factors by central units in
  $A$.
\end{definition}

\section{Motivation and Applications}

\begin{problem}
\label{problem2}
Let $A$ be a finite factorization domain and a $\KK$-algebra. Given $f\in A\setminus (U(A) \cup \{0\})$, compute all its factorizations of the form $f = c \cdot f_1 \cdots f_n$, where $c\in U(A)$ 
and $f_i \in A\setminus U(A)$ are irreducible.
\end{problem}

This paper is devoted in part to the algorithmic solution of Problem \ref{problem2} for a broad class of $G$-algebras. With this algorithm one can
approach a number of important problems, which we discuss in details.

Let $A$ be a $\KK$-algebra and $0 \neq L\subset A$ a finitely generated left ideal. 

\begin{problem}
\label{problem1}
Compute a proper left ideal $N\supsetneq L$.
\end{problem}
Unfortunately, it is not known in general, whether the problem of left maximality of a given ideal with respect to inclusion is decidable. Therefore we are interested in the local negative form of it. Namely, if Problem \ref{problem1} can be solved, then 
$L$ is not left maximal. Moreover, for any $N$
as above, we have a surjection from $A/L$ to its proper factor-module $A/N$,
in other words the exact sequence of left $A$-modules
\[
0 \to N/L \to A/L \to A/N \to 0
\]
which contributes to the knowledge of the structure of $A/L$.

Suppose that $f\in A\setminus\{0\}, f \notin L$ has finitely many
factorizations $f = g_i h_i$ up to multiplication by central units, where $i\in I, |I| < \infty$ and $g_i, h_i \notin U(A)$. We do not require that $g_i$ or $h_i$ are irreducible. Suppose, that $L \subsetneq L + Af \subsetneq A$.
Then $L + Af \subseteq \bigcap_{i\in I} (L + Ah_i)$, hence
there is a natural surjective homomorphism of left $A$-modules
\[
\frac{A}{L + Af} \to \frac{A}{\bigcap_{i\in I}( L + Ah_i )} \to 0
\]

{\bf Relations with solution spaces}: 
Let $\FF$ be an arbitrary (in particular, not necessarily finitely generated) left $A$-module, which one can think about as of the {\it space of solutions} for $A$-modules. Then, for fixed $\KK$ and an $A$-module $M$ one denotes
$\Sol(M, \FF) := \Hom_A(M, \FF)$, which is a $\KK$-vector space and 
an $\End_A(M)$-module \cite{WS}.

By invoking the Noether-Malgrange isomorphism \cite{WS}, we obtain for an
the natural injective map of $\KK$-vector spaces
\[
0 \to \Sol\left( \frac{A}{\bigcap_{i\in I} (L + Ah_i)} , \FF\right) \to \Sol\left(\frac{A}{L + Af}, \FF\right).
\]
The latter sheds light on the structure of the space of solutions of $A/(L+Af)$.

Note, that the following version of the left Chinese remainder theorem
for modules holds:
\begin{theorem}
\label{crt}
Let $A$ be a $\KK$-algebra, $I$ a finite set of indices and $\{L_i : i\in I\}$ are
 left ideals in $A$. Consider the homomorphism
\[
A/\bigcap_{i\in I} L_i \overset{\phi}{\longrightarrow} \bigoplus_{i\in I} A/ L_i, \quad a+\bigcap_{i\in I} L_i \mapsto (a+ L_1, \ldots, a+L_{|I|}).
\]
Then the following holds
\begin{itemize}
\item[1)] $\phi$ is injective
\item[2)] if $\forall i, j \in I, i\neq j$ holds $L_i + L_j = A$, then $\phi$ is surjective.
\end{itemize}
\end{theorem}
Of course, one can assume that $L_i$ are proper nonzero ideals.

In the second item of the Theorem \ref{crt} one says that the collection $\{L_i : i\in I\}$ is {\bf left comaximal}. Then $\phi$ is an isomorphism and one has a finite direct sum decomposition of the module $A/\bigcap_{i\in I} L_i$.
Hence, there is a 
direct sum decomposition of the solution spaces
\[
\Sol\left(A/\bigcap_{i\in I} L_i, \FF\right) = 
\Sol\left(\bigoplus_{i\in I} A/ L_i, \FF\right) = 
\bigoplus_{i\in I} \Sol\left(A/ L_i, \FF\right).
\]
Note, that the right hand side can be a direct sum even if the condition (2) is not satisfied, see Example \ref{exTsai}.

{\bf Another application}: Let $A$ be a domain and $S\subset A$ be multiplicatively closed Ore set. By Ore's Theorem the localization $S^{-1}A$
exists and there is an injective homomorphism $A \to S^{-1}A$ \cite{Bueso:2003}. 

A left ideal $L\subset A$ is called {\bf left $S$-closed} if $L^S = L$, where
$L^S := \{a\in A \mid \exists s\in S \; sa \in L \} \supseteq L$ is the {\bf $S$-closure of $L$}.
There is another characterization of $S$-closedness:
$L^S = \ker (A \to S^{-1}(A/L))$, where the latter homomorphism of $A$-modules is $a \mapsto 1^{-1} a+ S^{-1}L$.
Then $L^S/L$ is the $S$-torsion submodule of $A/L$ and $A/L^S$ has no $S$-torsion.

\begin{problem}
\label{problem3}
Given $S\subset A$ an Ore set and $L \subset A$, give an algorithm to compute $L^S$. 
\end{problem}

For a general $S$, it is unknown, whether $L^S$ is computable. If $A$ is the
$n$th Weyl algebra, $S=\KK[x_1,\ldots,x_n]\setminus\{0\}$ and $L\subset A$ 
has finite holonomic rank, then there is an algorithm \cite{TsaiWC, Tsaidiss} to compute $L^S$ (known as the {\bf Weyl closure of $L$}).

The factorization can be used in the process of computing $L^S$ as follows.
Let $A$ be an FFD. Given $\ell \in L$, one computes finitely many
factorizations $\ell = a_i b_i, i\in I$ for some finite indexing set $I$. Then let $J:= \{ j\in I \mid (a_j, b_j) \in S \times L \}$.
If $J \neq \emptyset$, one has $L + \{Ab_j : j \in J\} \subseteq L^S$. In such a way one obtains an approximation to $L^S$. Note, that $L^S = A$ if and only
if $L\cap S \neq \emptyset$.


\section{How to factor in {\fontsize{1em}{1em}{$G$}}-Algebras}


\subsection{General Algorithm}

In a recent publication \cite[Theorem 1.3]{bell2014noncommutative}, it was proven
that each $G$-algebra $\GG$ is a finite factorization domain.

In the same paper, an outline was given how one could find all
possible factorizations of an element in $\GG$. In this
section, we will provide a thorough description of an algorithm to
find all possible factorizations of an element in a  $G$-algebra $\GG$,
up to multiplication by central units.

For this, we need to make a further assumption on our field \KK, which
holds for most practical choices of \KK.\\

\noindent
\textbf{Assumption:} There exists an algorithm to determine if a
polynomial $p$ in $\KK[x]$ has roots in \KK. If $p$ has roots in \KK,
then this algorithm can produce all \KK-roots of $p$.\\

Recall, that $\{x^{\alpha} = x_1^{\alpha_1} \cdot \ldots \cdot x_n^{\alpha_n} : \alpha \in\NN_0^n \}$ is a $\KK$-basis of $\GG$. With respect to an admissible monomial ordering $\prec$ on $\GG$ we can uniquely write every $g \in \GG\setminus\{0\}$ as $g = c_{\alpha} x^{\alpha} + t_g$
with $c_{\alpha} \in \KK\setminus\{0\}$. Moreover, either $t_g=0$ or $x^{\beta} \prec x^{\alpha}$ for any summand
$c_{\beta} x^{\beta}$,  $c_{\beta} \neq 0$ of $t_g$. Then 
$\lm(g) = x^{\alpha}$ is the \emph{leading monomial} of $g$ and 
$\lc(g) = c_{\alpha}$ is the \emph{leading coefficient} of $g$. A polynomial $g\neq 0$
with $\lc(g)=1$ is called \emph{monic}.

It is important to recall \cite{LVDiss}, that $\lm(x^{\alpha} \cdot x^{\beta}) = x^{\alpha +\beta}$ holds $\forall \alpha,\beta \in\NN_0^n$ in a $G$-algebra.

\begin{algorithm}[t]
\caption{Factoring an element $g$ in a $G$-algebra \GG}
\label{alg:FactorGAlg}
\begin{flushleft}
\textit{Input:} $g \in \GG\setminus \KK$.\\
\textit{Output:} \(\{(g_1, \ldots, g_m) \mid  m\in \NN, g_i \in
\GG\setminus \KK \text{ for } i \in
\{ 1,\ldots , m\}, g_1\cdots g_m = g\}\) (up to
multiplication of each factor by a central unit).\\
\textit{Assumption:} An admissible monomial ordering $\prec$ on $\GG$ is fixed and $g$ is monic with respect to it.
\end{flushleft}
\begin{algorithmic}[1]
  \STATE $R := \{\}$
  \STATE \begin{align*}M := \{(p_1,\ldots,p_\nu) \mid &\nu \in \NN, p_i \in
  \{x_1,\ldots, x_n\}, \\
  &\lm(p_1\cdot\ldots\cdot p_\nu) = \lm(g)\}\end{align*}\label{ln:monomialCombinations}
  \FOR{$(p_1,\ldots, p_\nu) \in M$}
    \FOR{$i := 1$ \TO $\nu-1$}
      \STATE Set up an ansatz for the \KK-coefficients of  $a\cdot b = g$ with
      $\mathrm{lt}(a)=p_1\cdot\ldots\cdot p_i$ and $\mathrm{lt}(b) =
      p_{i+1}\cdot\ldots\cdot p_\nu$.\label{ln:ansatz}
      \STATE $F :=$ the reduced Gr\"obner basis w.r.t. an elimination ordering of
      the ideal generated by the coefficients of $a\cdot b-g$.
      \IF{$F\neq \{1\}$}
      \STATE $V :=$ Variety of $\langle F \rangle$ in an affine space over $\KK$. \label{ln:variety}
      \STATE $R := R \cup \{(a,b) \mid a, b \in \GG, a\cdot b =g$,
      where the coefficients of $a, b$ are given by $v\in V\}$\label{ln:tupleL}
      \ENDIF
    \ENDFOR
  \ENDFOR
  \IF{$R = \{\}$}
  \RETURN \{(g)\}
  \ELSE
  \STATE Recursively factor $a$ and $b$ for each $(a,b) \in R$.
  \ENDIF
  \RETURN R
\end{algorithmic}
\end{algorithm}

\begin{proof}[of Algorithm \ref{alg:FactorGAlg}]
Let us begin with discussing the termination. The set $M$ in line
\ref{ln:monomialCombinations} is finite, as it is a permutation of a
finite product of the variables in \GG. Since $\GG$ is a $G$-algebra, 
the set of total well-orderings on it, satisfying the Definition \ref{def: G-algebra}, 
is nonempty. By \cite{GomezRefiltering}, in this set there is a weighted degree total
ordering, say $\prec_w$ with strictly positive weights. Without loss of generality let us
assume this is the ordering we are working with. Thus for any monomial 
there are only finitely many monomials which are smaller with respect to $\prec_w$.
In particular, this applies to $\mathrm{lt}(a)$ and $\mathrm{lt}(b)$ in line \ref{ln:ansatz}. The
variety $V$ will be a finite set due to the fact that $\GG$ is an
FFD. Thus, the set $R$ in line \ref{ln:tupleL} will also be finite.
The recursive call will also terminate, since in each step we either
discover that we cannot refine our factorization any more, or we split
a given factor into two factors of strictly smaller degrees.

For the correctness discussion of our algorithm, we need to show that
we can calculate the variety $V$ in line \ref{ln:variety}. We know,
since $\GG$ is an FFD, that the ideal generated by $F$ is either zero-dimensional
over $\KK[x]$ or it is an intersection of such with a higher-dimensional ideal $H$, whereas 
the variety of $H$ does not contain points from an affine space over $\KK$. Hence we proceed
with the zero-dimensional component $F_0$ of $F$.

Our assumption above states that we can find all \KK-roots of a univariate polynomial. Since $F_0$  
is zero-dimensional, for any variable $x_i$ there is the corresponding univariate polynomial,
generating the principal ideal $F_0 \cap \KK[x_i]$. By backwards substitution, we obtain the entire \KK-variety of the ideal generated by $F_0$.
\end{proof}

\begin{example}
\label{ex1}
Let us consider the universal enveloping algebra $U(\mathrm{sl_2})$ of
$\mathrm{sl}_2$ \cite{dixmier1977enveloping},
represented by
\[\KK\langle e, f, h \mid fe=ef-h,he=eh+2e,hf=fh-2f \rangle.\]
In $U(\mathrm{sl_2})$, we want to factorize the element
\begin{align*}
p:=&e^3f+e^2f^2-e^3+e^2f+2ef^2-3e^2h-2efh-8e^2\\
&+ef+f^2-4eh-2fh-7e+f-h.
\end{align*}
We fix the lexicographic ordering on $U(\mathrm{sl_2})$, i.e. the
leading term of $p$ is $e^3f$.

Therefore the set $M$ in line \ref{ln:monomialCombinations} is given as
\[M := \{(e,e,e,f),(e,e,f,e),(e,f,e,e),(f,e,e,e)\}.\]

When choosing $(e,e,e,f)$, for $i=1$ one obtains the factorization
$$p = (e+1)\cdot(e^2f + ef^2-3eh - 2fh -e^2+f^2-7e+f-h).$$

By picking $(e,e,f,e)$, for $i= 3$ one obtains two more
factorizations, namely
$$p=(e^2f + 2ef -2eh -e^2 -4e + f -2h -3)\cdot(e+ f)$$
and
$$p=(e^2f + ef^2 -2eh -e^2 +f^2 -3e -f -2h)\cdot(e + 1).$$

All the other combinations either produce the same factorizations or
none.

When recursively calling the algorithm for each factor in the found
factorizations, we discover that the first two factorizations have a
reducible factor. In the end, one obtains the following two distinct
factorizations of $p$ into irreducible factors:
\begin{align*}
  p =& (e^2f + ef^2 -2eh -e^2 +f^2 -3e -f -2h)\cdot(e + 1)\\
    =& (e+1)\cdot (ef-e+f-2h-3)\cdot (e+f).
\end{align*}
\end{example}

\subsection{Implementation}

We have developed an experimental implementation of Algorithm \ref{alg:FactorGAlg} in the computer
algebra system {\sc Singular} \cite{DGPS}. We will make it available as part of {\tt
  ncfactor.lib}
. Our newly implemented procedures factorize elements in any
$G$-algebra, whose ground field is $\mathbb{F}(q_1,\ldots,q_n)$, where
$\mathbb{F}$ is either $\QQ$ or a finite prime field and $q_i$ are transcendental over $\mathbb{F}$.

We designed the software in a modular way, so that during runtime our
function checks if a more efficient factorization algorithm is
available for the specific given $G$-algebra and/or input polynomial. If
this is the case, the input is re-directed to this function. In this
way, the user can call the general function to factor elements in any one of the
supported $G$-algebras, and runs the available optimized algorithms,
where available, without calling them individually.

\subsection{Possible Improvements}

Algorithm \ref{alg:FactorGAlg} solves the problem of finding all
possible factorizations of an element in a $G$-algebra, but it will not be very efficient in general. This is not only due to the
complexity of the
necessary calculation of a Gr\"obner basis \cite{mayr1982complexity},
but also the size of the set $M$ is a bottleneck. In
\cite{giesbrecht2014factoring,Giesbrecht2015}, an algorithm for
factoring elements in the $n$th Weyl algebra is presented, which
is similar to Algorithm \ref{alg:FactorGAlg}. The main difference is
that the
$\ZZ^n$-graded structure is utilized. There, the homogeneous
polynomials of degree zero form a $\KK$-algebra $A_n^{(0)}$, which is isomorphic to a commutative multivariate polynomial ring. The set of homogeneous polynomials of degree $z \in \ZZ^n \setminus \{0^n\}$ has the structure of a cyclic $A_n^{(0)}$-bimodule. Hence, factorization of homogeneous polynomials with respect to the $\ZZ^n$-grading reduces to factoring commutative polynomials with minor additional combinatorial steps. 
An inhomogeneous polynomial $f$ has now the highest graded part $\alpha(f)$ 
and the lowest graded part $\omega(f)$, both of them rather polynomials than 
monomials. Hence $\alpha(f), \omega(f)$ have potentially smaller numbers of different factorizations than the permutations of the leading term collected in $M$ in Algorithm \ref{alg:FactorGAlg}. Indeed, it suffices to consider firstly factorizations into two polynomials and for each candidate pair an 
 ansatz is made for the graded terms between the highest and the
 lowest graded parts. This means, that the set $M$ has smaller size in
 general when using this technique. Additionally, this approach takes
 the lowest graded part into account, which allows to eliminate
 certain invalid cases beforehand.
The performance increase is reflected by the benchmarks presented in
\cite{giesbrecht2014factoring,heinle2013factorization}.

Hence, for practical implementations
of Algorithm \ref{alg:FactorGAlg},
one should examine each possible $G$-algebra separately and take
advantage of potential extra structure, like the presence of nontrivial 
$\ZZ^n$-grading or an isomorphism to an algebra with this structure.

We will conclude this section by summarizing the conditions that can
lead to an improved version of Algorithm \ref{alg:FactorGAlg}. Let $A$
be a $\KK$-algebra, which possesses a nontrivial (i.e. not all weight vectors are zero) $\ZZ^n$-(multi)grading. Then one can infer the following additional information:
\begin{enumerate}
\item For $z \in \ZZ^n$, $A_z:=\{ a\in A: \deg(a)=z\}\cup \{0\}$ is a $\KK$-vector space. Moreover, $\oplus_z A_z = A$ and $A_i A_j \subseteq A_{i+j}$ for all $i, j \in \ZZ^n$.
\item $A_{0^n}$, the graded part of degree zero, is a $\KK$-algebra itself (since $A_0 A_0 \subseteq A_0$).
\item For $z \in \ZZ^n \setminus \{0^n\}$, the $z$-th graded part $A_z$ is an $A_{0^n}$-bimodule (since $A_0 A_z, A_z A_0 \subseteq A_z$).
\end{enumerate}

In order to be useful for factorizing purpose, this grading should have the following properties:
\begin{enumerate}
\setcounter{enumi}{3}
\item The graded part of degree zero,  $A_{0^n}$, which is a $\KK$-algebra,
is additionally an FFD with "easy" factorization, preferably the
commutative polynomial ring. Furthermore, for keeping the set $M$ in
Algorithm \ref{alg:FactorGAlg} small, it would be desirable if in $A_{0^n}$
a randomly chosen polynomial is irreducible with high probability.
\item The irreducible elements in $A_{0^n}$, that are reducible in
  $A$, can be identified and factorized in an efficient manner. Preferably, one has a finite number of monic elements of such type.
\item For $z \in \ZZ^n \setminus \{0^n\}$, the $z$-th graded part
  $A_z$ is a finitely generated $A_{0^n}$-bimodule, preferably a
  cyclic bimodule.
\end{enumerate}
Then the Algorithm \ref{alg:FactorGAlg} can be modified along the lines 
of algorithms from \cite{giesbrecht2014factoring,Giesbrecht2015}, which we
have also sketched above. Let us illustrate this approach by a concrete example.

\begin{example}
\label{ex2}
As in Example \ref{ex1}, let $A=U(\mathrm{sl_2})$, that is
\[ A = \KK\langle e, f, h \mid fe=ef-h,he=eh+2e,hf=fh-2f \rangle.\]
Af first, let us determine which gradings are possible. Let $w_e, w_f$
and $w_h$ be the weights of the variables, not all zero. The two last relations of $A$ imply that $w_h=0$, and the first one implies $w_e + w_f = w_h = 0$, that is $w_f = - w_e$. Hence a $\ZZ$-grading $(w_e, w_f, w_h) = (1,-1,0)$ is enough for our purposes, since $A_0 = \KK[ef,h]$ is commutative and the $z$-th graded part is a cyclic $A_0$-bimodule, generated by $e^z$ if $z>0$ and by
$f^{|z|}$ otherwise. This property guarantees, that $\forall r \in \KK[ef,h]$
and $\forall z\in \NN$ there exists $q_1, q_2 \in \KK[ef,h]$, such that 
$r e^z = e^z q_1$ and $e^z r = q_2 e^z$ and the same holds for the multiplication by $f^z$. Note, that $\deg(q_i)=\deg(r)$.

We claim that the only monic irreducible elements in $A_0$, which are
reducible in $A$, are given by $ef$ and $ef-h$. The proof to this
claim is similar to the one for \cite[Lemma 2.4]{Giesbrecht2015},
which we outline here: Let $p$ be an irreducible element in $A_0$,
which reduces into $p=\varphi \cdot \psi$ in $A$, where $\varphi,
\psi\in A\setminus \KK$ are monic. Since $A$ is a domain, the factors  $\varphi,\psi$ are
homogeneous with $\deg(\varphi) = k$ and $\deg(\psi)=-k$ for some $k
\in \ZZ$. If $|k|>1$ or $k =0$, $p$ would be reducible in $A_0$, which violates
our assumption. Hence only $k=1$ is possible. If any of $\varphi$ or
$\psi$ would have a non-trivial $A_0$ factor, we would obtain again
that $p$ is reducible in $A_0$. This leaves as only options $p=ef$
or $p= fe = ef-h$, as claimed. Thus, we have shown that irreducible
elements in $A_0$, which are reducible in $A$, can be easily
identified and factored.

Now consider the same polynomial $p$ as in Example \ref{ex1}. With respect to the $(1,-1,0)$-grading it decomposes into the following graded parts:
$\alpha(p) = -e^3$, $\omega(p) = f^2$ (as we see, in this case we have monomials in graded parts, while in general rather polynomials appear) and the intermediate parts are
\[
\underbrace{e^3 f-3e^2 h-8e^2}_{\deg: 2} +
\underbrace{e^2 f-4eh-7e}_{\deg: 1} +
\]
\[
+\underbrace{e^2 f^2-2efh+ef-h}_{\deg: 0} +
\underbrace{2ef^2-2fh+f}_{\deg: -1}.
\]
Among the factorizations of $\alpha(p) = -e^3$ and $\omega(p) = f^2$ into two factors, consider the case $(-e^2) \cdot e$ and $f \cdot f$. Thus, we're looking for $a, b \in A$ with $\alpha(a)=e^2, \omega(a)=f$ and  
$\alpha(b)=e, \omega(b)=f$ and $p=ab$ holds. In $b$ we have only one possible intermediate graded part $b_0(ef,h)$, namely of degree 0 since $\deg\alpha(b) = 1$ and $\deg \omega(b)=-1$. In $a$ we have to specify the parts of degrees 1 resp. 0, that is $a_1(ef,h)\cdot e$ resp. $a_0(ef,h)$. After the multiplication, we obtain the following graded decomposition of intermediate graded terms of $ab$:
\[
\underbrace{-e^2 b_0 + a_1 e^2}_{\deg: 2} + 
\underbrace{a_1 e b_0 + a_0e - e^2 f}_{\deg: 1} + 
\]
\[
+
\underbrace{a_1 ef + a_0 b_0 + ef -h}_{\deg: 0} +
\underbrace{f b_0 + a_0 f}_{\deg: -1}.
\]
By fixing the maximal possible degree of $a_0, a_1, b_0 \in
\KK[ef,h]$, we can create and solve a system of equations which the
coefficients of $a_0, a_1, b_0$ have to satisfy. In this example an ansatz
in terms of $1, h, ef$, i.e. 9 unknown coefficients, leads to the system
of 18 at most quadratic equations, which leads to the unique solution: $b_0(ef,h)=0$, 
 $a_0 (ef,h)=2ef-2h-3$ and $a_1 (ef,h) = ef-h-2$.  Substituting the polynomials, we arrive at the following factorization with polynomials sorted
 according to the grading:
  \[
 p = (-e^2+e^2f-2eh-4e+2ef-2h-3+f) \cdot (e+f)
 \]
This is already known to us from the Example \ref{ex1}.
In an analogous way one can address other factorizations. Note, that in the ansatz we made, significantly less variables for unknown coefficients
and a system of less equations of smaller total degree were used, compared to the general Algorithm. 
\end{example}

\section{The Factorized Gr\"obner Basis Algorithm for
  {\fontsize{1em}{1em}{$G$}}-Algebras}

In what follows, by the term ideal we always mean a left ideal (unless
otherwise specified).

The factorized Gr\"obner approach has been studied extensively for the
commutative case
\cite{czapor1989solving,czapor1989solvingmultivariate,davenport1987looking,
grabe1995factorized, grabe1995triangular}, and implementations are e.g.
provided in the computer algebra systems \textsc{Singular} \cite{DGPS} and
\textsc{Reduce} \cite{hearnreduce}.

The leading motivation is to split a Gr\"obner basis computation into
smaller pieces with respect to the degrees of their generators. The union
of the varieties of the ideals generated by these smaller pieces
equals the variety of the initial system.

In the commutative case, there is also a way to constrain the solution
space. One can provide an extra set of elements, that should not be
reducible by the computed Gr\"obner basises. In this way, one excludes
certain unwanted solutions, which is useful in practice.

The search for varieties in the commutative case translates to the
search for solutions in the non-commutative case: All $G$-algebras are finite factorization
domains and a general factorization algorithm via
Algorithm \ref{alg:FactorGAlg} is given. Many of them are abstractions
of algebras of operators, and one is interested to find common
solutions of certain sets of operators, written as polynomials. Right hand factors of
elements correspond to partial solutions, and hence a split similar to
the commutative case is helpful to recover partial solutions. Motivated by this observation,
we attempt to generalize the factorized Gr\"obner basis algorithm to
the $G$-algebra case in this section. Our algorithm includes the possibility to
introduce constraints, similar to the methods in the commutative case.  

Unfortunately, not all nice properties transfer into the
non-commutative case, as the following example depicts.

\begin{example}
\label{exTsai}
In the commutative case, one has the property that the radical
of the input ideal will be equal to the intersection of the radicals
of all ideals computed by the factorized Gr\"obner basis algorithm.

We will show via a counter-example that we do not have the same
property for $G$-algebras.

Consider
\begin{align*}
p =& (x^6+2x^4-3x^2)\partial^2-(4x^5-4x^4-12x^2-12x)\partial\\
& +
(6x^4-12x^3-6x^2-24x-12) \in A_1.
\end{align*}
This polynomial appears in \cite{TsaiWC} and has two different factorizations, namely
\begin{align*}
p =
&(x^4\partial-x^3\partial-3x^3+3x^2\partial+6x^2-3x\partial-3x+12)\cdot
\\
&(x^2\partial+x\partial-3x-1)\\
 =&
 (x^4\partial+x^3\partial-4x^3+3x^2\partial-3x^2+3x\partial-6x-3)\cdot
 \\
& (x^2\partial-x\partial-2x+4)
\end{align*}
A reduced Gr\"obner basis of $\langle
x^2\partial+x\partial-3x-1\rangle \cap \langle x^2\partial-x\partial-2x+4
\rangle$, computed in \textsc{Singular} \cite{DGPS}, is given by
\begin{align*}
\{&3x^5\partial^2+2x^4\partial^3-x^4\partial^2-12x^4\partial+x^3\partial^2-2x^2\partial^3+16x^3\partial\\
& +9x
^2\partial^2+18x^3+4x^2\partial+4x\partial^2-42x^2-4x\partial-12x-12,\\
&2x^4\partial^4-2x^4\partial^3+11x^4\partial^2+12x^3\partial^3-2x^2\partial^4-2x^3\partial^2\\
&+10x^
2\partial^3-44x^3\partial-17x^2\partial^2+64x^2\partial+12x\partial^2+66x^2\\
&+52x\partial+4\partial^2-168x-16
\partial-60\}.
\end{align*}
\end{example}

Hence, one main difference will be that we do not claim that the
union of all solutions of our smaller pieces in the factorizing
Gr\"obner basis algorithm will always be equal to all common solutions of
the initial set of polynomials. In general, we will only find a subset
of all solutions using our method. However, in this example, the
space of holomorphic solutions of the differential equation associated to $p$ in fact
coincides with the union of the solution spaces of the two generators
of the intersection stated above.

\begin{definition}
  \label{def:constrainedGroebnerTuple}
  Let $B,C$ be finite subsets in $\GG$. We call the tuple $(B,C)$ a
  \textbf{constrained Gr\"obner tuple}, if $B$ is a Gr\"obner
  basis of $\langle B \rangle$, and  $\mathrm{NF}(g,B)\neq 0$ for
  every $g \in C$. We call a constrained Gr\"obner tuple \textbf{factorized},
  if every $f \in B$ is either irreducible or has a unique irreducible
  left divisor.
\end{definition}

It is possible to strengthen the assumptions on a factorized
constrained Gr\"obner tuples by only allowing completely irreducible
elements in $B$, which might be preferable depending on the concrete
problem. However, in our application, we allow elements with
only one factorization. In this way, we increase the number of
solutions we can find for a certain system $B \subset \GG$ by using
our generalized factorized Gr\"obner basis algorithm. This methodology
also appears in the context of semifirs, where the concept of so called block
factorizations or cleavages has been introduced to study the
reducibility of a principal ideal \cite[Chapter 3.5]{cohn2006free}.

\begin{algorithm}[t]
\caption{Factorized Gr\"obner bases Algorithm for $G$-Algebras (FGBG)}
\label{alg:FacGSTD}
\begin{flushleft}
\textit{Input:} $B := \{f_1, \ldots, f_k\} \subset \GG$, $C :=
\{g_1,\ldots, g_l\} \subset \GG$.\\
\textit{Output:} $R := \{(\tilde B, \tilde C) \mid (\tilde B, \tilde
C)\text{ is factorized constrained Gr\"obner tuple}\}$ with
$\langle B \rangle \subseteq \bigcap_{(\tilde B, \tilde C) \in
  R}\langle \tilde B\rangle$ \\
\textit{Assumption:} All elements in $B$ and $C$ are monic.
\end{flushleft}
\begin{algorithmic}[1]
\FOR{$i= 1$ \TO k}\label{ln:beginAlterInitialGeneratorSet}
\IF{$f_i$ is reducible}
\STATE $M := \{(f_i^{(1)}, f_i^{(2)} \mid 
f_i^{(1)}, f_i^{(2)} \in \GG\setminus \KK, \mathrm{lc}(f_i^{(1)})= \mathrm{lc}(f_i^{(2)}) = 1,
 f_i^{(1)}\cdot f_i^{(2)}  = f_i, f_i^{(1)} \text{ is
   irreducible}\}$\label{ln:gb_element_factor_list}
\IF{there exists $(a,b),(\tilde a, \tilde b) \in M$ with $\tilde
a \neq a$}\label{ln:existenceCheckDifferentLeftDiv}
\RETURN $$\bigcup_{(a,b) \in M} \mathrm{FGBG}\left((B\setminus \{f_i\})\cup \{b\}, C \cup
\bigcup_{(\tilde a, \tilde b) \in M\atop b \neq \tilde b} \{\tilde b\} \right)$$\label{ln:recCall1}
\ENDIF
\ENDIF
\ENDFOR\label{ln:endAlterInitialGeneratorSet}
\STATE $ P := \{(f_i,f_j)\mid i,j \in \{1,\ldots, k\}, i < j\}$\label{ln:baseCaseBuchBerger}
\WHILE{ $P \neq \emptyset$} \label{ln:BuchbergerStart}
\STATE Pick $(f,g) \in P$
\STATE $P := P \setminus \{(f,g)\}$
\STATE $s:=$ S-polynomial of $f$  and $g$
\STATE $h := \mathrm{NF}(s, B)$
\IF{$h \neq 0$}
\IF{$h$ is reducible}
\RETURN $\mathrm{FGBG}(B \cup \{h\}, C)$
\ENDIF
\STATE $P := P \cup \{(h,f) \mid f \in B\}$
\STATE $B:= B \cup \{h\}$
\ENDIF
\IF{there exists $i \in \{1,\ldots, l\}$ with $\mathrm{NF}(g_i,B) = 0$}
\RETURN $\emptyset$
\ENDIF
\ENDWHILE \label{ln:BuchbergerEnd}
\RETURN $\{(B,C)\}$
\end{algorithmic}
\end{algorithm}

\begin{proof}[of Algorithm \ref{alg:FacGSTD}]
We will first discuss the termination aspect of Algorithm
\ref{alg:FacGSTD}. Since $M$ as calculated in line
\ref{ln:gb_element_factor_list} is of finite cardinality, the
existence check in line \ref{ln:existenceCheckDifferentLeftDiv}
can be done in a finite number of steps. Line \ref{ln:recCall1}
consists of a finite number of recursive calls to FGBG. The algorithm reaches this
line if there is an element $f$ in $B$, which is reducible and has a 
non-unique irreducible left divisor. In each recursive call, the
algorithm is called with an altered version of the set $B$, where
$f$ is being replaced in $B$ by $b \in \GG$, where $b$ is chosen
such that there exists an irreducible $a$ in $\GG$ with $f =
ab$. Therefore, after a finite depth of recursion, FGBG will be
called with a set $B$ containing elements that are either irreducible or
have an unique irreducible left divisor. We can make this assumption
on $B$ when FGBG reaches line \ref{ln:baseCaseBuchBerger}. Lines
\ref{ln:BuchbergerStart}--\ref{ln:BuchbergerEnd} describe the
Buchberger algorithm to compute a Gr\"obner basis, with two
differences:
\begin{enumerate}
\item If the normal form $h$ of an S-polynomial with respect to $B$ is not
  0, we check $h$ for reducibility. If $h$ is reducible, we call FGBG
  recursively, adding $h$ to $B$.
\item We check the system for consistency, i.e. if there is an element
  in $C$ that reduces with respect to $B$, we return the empty set.
\end{enumerate}
Each recursive call will terminate, since we add an element to $B$
that will reduce an S-polynomial to zero, which could not be reduced
to zero before.

For the correctness discussion, one observes that lines
\ref{ln:beginAlterInitialGeneratorSet}--\ref{ln:endAlterInitialGeneratorSet}
serve the purpose to split the computation based on the reducibility
of the elements in the initial set $B$. If an element $f \in B$
factorizes in more than one way, we recursively call FGBG with
$(B\setminus \{f\}) \cup
\{b\}$
 as the generator set for each
maximal right hand factor $b$ of $f$. Hence, the left
ideal generated by $(B\setminus\{f\})\cup \{b\}$ will contain $\langle
B\rangle$,
and thus $\langle B \rangle$ will be contained in the intersection of
all of them, as
required.

As already mentioned in the termination discussion, lines
\ref{ln:BuchbergerStart}--\ref{ln:BuchbergerEnd} describe the
Buchberger algorithm. After computing an S-polynomial $h$, we check
for its reducibility. If there is more than one maximal right factor
$r$ of $h$, we call FGBG recursively and add $h$ to our set $B$. Here,
we have again a guarantee that the left ideal generated by $B$ is a
subset of the left ideal generated by $B\cup\{h\}$.

The additional constraints that we impose on each recursive call
enable us to minimize our computations, but do not violate the subset
property. In the end, it is ensured that in all computed constrained
Gr\"obner tuples $(\tilde B, \tilde C)$, no element in $C$ lies in the
left ideal generated by $\tilde B$.
\end{proof}


\begin{example}
  Let us execute FGBG on an example. Let
  \begin{align*}
    B := \{&\partial^4+x\partial^2-2\partial^3-2x\partial+\partial^2+x+2\partial-2,\\
         &
    x\partial^3+x^2\partial-x\partial^2+\partial^3-x^2+x\partial-2\partial^2-x+1\}
  \end{align*}
  be a subset of the first Weyl algebra $A_1$. We assume that
  $C:=\{\partial -1\}$, and that our ordering is the degree
  reverse lexicographic one with $\partial > x$. This example is taken
  from the \textsc{Singular} manual \cite{DGPS} (and it is a Gr\"obner
  basis for the left
  ideal $\langle \partial^2+x  \rangle\cap \langle \partial -1
  \rangle$; hence we would expect the output with our chosen $C$ to be $\langle \partial^2+x  \rangle$). Each element factors
  separately as
  \begin{align*}
    f_1 := &\partial^4+x\partial^2-2\partial^3-2x\partial+\partial^2+x+2\partial-2\\
    =& (\partial^3+x\partial-\partial^2-x+2)\cdot
    (\partial-1)\\
    =& (\partial-1)\cdot(\partial^3+x\partial-\partial^2-x+1)
  \end{align*}
  respectively
  \begin{align*}
    f_2 :=& x\partial^3+x^2\partial-x\partial^2+\partial^3-x^2+x\partial-2\partial^2-x+1\\
    =& ( x\partial^2+x^2+\partial^2+x-\partial-1)\cdot(\partial-1)\\
    =& (x\partial-x+\partial-2)\cdot (\partial^2+x).
  \end{align*}
  Hence, in line \ref{ln:recCall1}, FGBG will return two recursive calls
  of itself, namely
  \begin{itemize}
  \item $\mathrm{FGBG}(\{\partial -1,
    f_2\},\{\partial-1,\partial^3+x\partial-\partial^2-x+1
    \})$
  \item $\mathrm{FGBG}(\{\partial^3+x\partial-\partial^2-x+1,f_2\},C)$
  \end{itemize}
  For simplicity, we will ignore the first call, as $C$ contains
  $\partial -1$, which also appears in the generator list.

  Furthermore, the new element
  $b_1 :=\partial^3+x\partial-\partial^2-x+1$ only has one
  possible factorization. Therefore, we consider now the
  factorizations of $f_2$. This leads again in line \ref{ln:recCall1}
  to two recursive calls:
  \begin{itemize}
    \item $\mathrm{FGBG}(\{b_1, \partial -1\}, \{\partial -1,
      \partial^2+x\})$
    \item $\mathrm{FGBG}(\{b_1, \partial^2+x\},C)$
  \end{itemize}
  As before, we can ignore the first recursive call. Thus, we are
  left with $(\{b_1, \partial^2+x\},C)$ to proceed on line
  \ref{ln:baseCaseBuchBerger}. 

  The normal form of the S-polynomial of $b_1$ and $\partial^2 +x$ is
  equal to zero. Further, the normal form of $b_1$ with respect to
  $\langle \partial^2+x\rangle$, is equal to zero, i.e. $\partial^2+x$ is a right
  divisor of $b_1$. Hence, we can omit $b_1$ and our complete Gr\"obner basis is given by $\{\partial^2+x\}$. 
  Since $\mathrm{NF}(\partial -1, \langle\partial^2+x\rangle) \neq 0$, our algorithm returns
  $\{(\{\partial^2+x\},C)\}$ as final output.

  Note, that if we would have chosen $C = \emptyset$ in the beginning,
  the output of our algorithm would have been $$\{(\{\partial
  -1\},\{b_1\}),(\{\partial^2 + x\},\{\partial -1\})\},$$
i.e. we recover $\langle B \rangle =\langle \partial^2+x  \rangle\cap \langle \partial -1
  \rangle$ in this case.


  

\end{example}

\begin{remark}
  One can also insert an early termination criterion inside Algorithm
  \ref{alg:FacGSTD}, namely after at least one factorized constrained
  Gr\"obner tuple has been found. This is in the commutative case
  motivated by the fact that in practice users are often not
  interested in all the elements in a variety but would be content with at least
  one. For example, the
  computer algebra system \textsc{Reduce} can be
  instructed to stop after finding one factorized Gr\"obner basis (see
  \cite{hearnreduce}). In the non-commutative case, we can only hope
  for partial solutions in general, but a mechanism to stop a
  computation once at least one is found is also desirable.
\end{remark}

\section{Conclusions}

An algorithm for factoring elements in $G$-algebras, where the
underlying field $\KK$ has the property that we are able to extract
all possible $\KK$-roots of any polynomial in $\KK[x]$, has been
shown (Algorithm \ref{alg:FactorGAlg}).

This algorithm and the FFD-property of $G$-algebras enable us to
propose a generalization of the factorized Gr\"obner basis algorithm
for $G$-algebras (Algorithm \ref{alg:FacGSTD}).

A future work would be to identify improvements to Algorithm
\ref{alg:FactorGAlg} for practically interesting $G$-algebras.
This has been studied e.g. for partial $q$-differential, differential and
difference operators in \cite{giesbrecht2014factoring,Giesbrecht2015},
where the $\ZZ^n$ graded structure resp. a certain embedding has been
utilized.
In the meantime, we have implemented the unimproved version in the
{\sc Singular} library {\tt ncfactor.lib}, which will be made
available shortly. Our implementation identifies beforehand if an
improved method is already included in {\tt ncfactor.lib} for a specific algebra and, if
this is the case,
re-directs the input there. This modular design allows us to update
the function once an improved algorithm is available for a certain
$G$-algebra. The use of the function stays the same after such an update.

Another interesting future direction would be to characterize further the
connection between the solution space of a polynomial system $B
\subset \GG$ and  the union of the solution spaces of
the output of Algorithm \ref{alg:FacGSTD} when called with
$B$. Especially, it would be interesting to identify properties of
$\GG$ and $B$, under which both spaces coincide.

An implementation of Algorithm \ref{alg:FacGSTD} would also be of
practical interest, which the authors intend to provide in the near
future.

\section{Acknowledgments}

We thank Eugene Zima, Jason P. Bell and Mark Giesbrecht for insightful discussions we had with them. We are grateful to Wolfram Koepf for
the information on internals of {\sc Reduce}.

%
\bibliographystyle{abbrv}
\bibliography{fac_g_std}  
\end{document}